\documentclass[ECP]{ejpecp}
\newcommand{\bbn}{\mathbb N}
\newcommand{\bbe}{\mathbb E}

\newcommand{\bbp}{\mathbb P}

\newcommand{\cala}{\cal A}
\SHORTTITLE{Longest Common Subsequences in Random Words}
\TITLE{Closeness to the Diagonal for Longest Common Subsequences in Random Words\thanks{Partially supported by the 
grant \# 246283 from the Simons Foundation and by a Simons Fellowship grant \# 267336}}
\AUTHORS{Christian Houdr\'e$^\dag$ \and Heinrich Matzinger\thanks{Georgia Institute of Technology, 
School of Mathematics, Atlanta, Georgia, 30332-0160}}
\AMSSUBJ{05A05 ; 60C05 ; 60F10}

\SUBMITTED{April 2016}
\ACCEPTED{April 20, 2016}
\VOLUME{0}
\YEAR{2016}
\PAPERNUM{4029}
\DOI{vVOL-PID}

\ABSTRACT{The nature of 
the alignment with gaps corresponding to a longest 
common subsequence (LCS) of two independent iid random sequences drawn from 
a finite alphabet is investigated.  It is shown that such 
an optimal alignment typically matches pieces of 
similar short-length.   This is of importance in understanding  
the structure of optimal alignments of two sequences.  
Moreover, it is also shown that 
any property, common to two subsequences, 
typically holds in most parts of the optimal 
alignment whenever this same property 
holds, with high probability, for strings of similar short-length.  
Our results should, in particular, prove useful for simulations since 
they imply that the re-scaled two dimensional representation of a LCS gets uniformly 
close to the diagonal as the length of the sequences grows without bound.}

\begin{document}
\section{Introduction}

Let $x$ and $y$ be two finite strings.  A common subsequence
of $x$ and $y$ is a subsequence of both $x$ and $y$, while a 
longest common subsequence (LCS) 
of $x$ and $y$ is a common subsequence of maximal length.

It is well known that common subsequences can be represented 
via \textit{alignments with gaps} as illustrated, next, on some examples:  
First take the binary strings $x=0010$ and $y=0110$. 
A common subsequence is $01$, which can be represented as an 
alignment with gaps as follows:  
The common letters are aligned together, while 
each letter not appearing in the common subsequence 
is aligned with a gap. Several alignments 
can thus represent the same common subsequence and in this first example, 
an alignment corresponding to the common subsequence
$01$ is given by
\begin{equation}
\label{bad}
\begin{array}{c|c|c|c|c|c|c|c}
x& &0&0&1&0& & \\\hline
y& & &0&1& &1&0
\end{array}
\end{equation}
another one is given by 
\begin{equation}
\label{bad2}
\begin{array}{c|c|c|c|c|c|c|c|c}
x& &0&0&1&0& & & \\\hline
y& &0& &1& &1&0
\end{array}
\end{equation}

Above, the LCS is not $01$ but rather 
$010$.  An alignment corresponding to a LCS
is called an \textit{optimal alignment} (OA) 
or is said to be \textit{optimal}.  
Neither \eqref{bad} nor \eqref{bad2} represent optimal alignments, 
but an optimal alignment is given by:
\[
\begin{array}{c|c|c|c|c|c|c}
x& &0&0&1& &0 \\\hline
y& & &0&1&1&0
\end{array}
\]
which, again, is clearly not unique.  
Here the LCS of $x$ and $y$
is $LCS(x;y)=010$, which has length three, a fact denoted  
by $|LCS(x;y)|=3$.  Here is another example: 
let $x=christian$ and $y=krystyaan$.  Then, 
$LCS(x;y)=rstan$ and an alignment with gaps representing the LCS is:
\begin{equation}
\label{alignmenthein}
\begin{array}{c|c|c|c|c|c|c|c|c|c|c|c|c|c|c|c|c}
&c& &h&r&i& &s&t&i& &a& &n \\\hline
& &k& &r& &y&s&t& &y&a&a&n 
\end{array}
\end{equation}
Again, all the letters which are part of the LCS are aligned with each 
other, while the other letters are aligned with gaps.  
In the above alignment of $x$ and $y$,   
$x_5x_6x_7x_8x_9=stian$ is aligned (with gaps) with 
$y_4y_5y_6y_7y_8y_9=styaan$, and we say that the 
integer interval $[5,9]_{\bbn}$ is aligned 
with $[4,9]_{\bbn}$. 
Alternatively, we say that $[5,9]_{\bbn}$ gets mapped 
to $[4,9]_{\bbn}$ by the alignment we consider, meaning that 
the following two conditions are satisfied:
\begin{itemize}
\item[i)] The letters $x_5x_6x_7x_8x_9$ 
are all aligned exclusively with gaps or with letters
from the string $y_4y_5y_6y_7y_8y_9$.
\item[ii)] The letters from $y_4y_5y_6y_7y_8y_9$
 are all aligned with gaps or with 
letters from the substring $x_5x_6x_7x_8x_9$.\end{itemize}

To emphasize our terminology, we see that in the 
alignment \eqref{alignmenthein}, 
$[1,4]_{\bbn}$ is aligned 
with $[1,2]_{\bbn}$ ($[1,4]_{\bbn}$ is also aligned with $[1,3]_{\bbn}$).  
In other words, in an alignment a piece of $x$ gets aligned with a piece
of $y$ if and only if the letters 
from the piece of $x$ which get aligned to letters
get only aligned with letters from the piece of $y$ and vice versa.  
Longest common subsequences and optimal alignments 
are important tools used in Computational Biology and 
Computational Linguistics for strings matching, e.g., 
see \cite{C}, \cite{RRS}, \cite{SK}, 
and \cite{W2}.

\textit{
Throughout this paper, 
$X=X_1\cdots X_n$ and $Y=Y_1\cdots Y_n$ are   
two independent random strings/words where  
$(X_i)_{i\ge 1}$ and $(Y_i)_{i\ge 1}$ are two independent 
iid sequences drawn from a finite alphabet
$\cala$. (No assumption, besides its non triviality, is made on 
$X_1$.)  Further, and again throughout, $LC_n$ denotes the 
length of the LCSs of $X$ and $Y$, i.e., 
$LC_n:=|LCS(X_1X_2\cdots X_n; Y_1Y_2\cdots Y_n)|$.}

To further put our framework in context, also note that 
$2(n-LC_n)$ is a version of the edit/Levenshtein distance used in computer science.  
It is equal to the minimal number of insertions and deletions to change 
either string/word  into the other.

A well known result of Chv\'atal and Sankoff \cite{CS} asserts that, when scaled by 
$n$, $\bbe LC_n$ converges to a constant $\gamma^*\in (0,1)$ 
(which depends on the 
size of the alphabet and on the law of $X_1$) given, via superadditivity, by 
$\gamma^*=\sup_{n\ge 1} \bbe LC_n/n$.  However, to this day, even 
in the uniform binary case, 
the exact value of $\gamma^*$ is unknown.  Moreover, Alexander~\cite{A}  
determined the rate of convergence to $\gamma^*$ showing 
that there exists an absolute constant $K_A$, independent of 
$n$, of the size of the finite alphabet $\cala$, (and of the law of $X_1$), 
such that for all $n\ge 1$,
\begin{equation}\label{alex}
\gamma^* n-K_A\sqrt{n\log n} \le \bbe LC_n 
\le\gamma^* n.  
\end{equation}

This law of large numbers (with rate) gives the first order behavior of $LC_n$, and 
the next problem of interest is to study the order of the variance 
(or more generally, of the centered absolute moments) of $LC_n$.  
This order has been found to be linear in the length of the sequences 
in various instances (e.g., \cite{HLM}, \cite{HM}, \cite{LM}, \cite{HMa}). 
This is, in particular, the case for iid binary sequences with 
zeros and ones having very different probabilities or for some classes of models 
``as close as one wants" to the uniform iid one \cite{AHM2}.  
In all the known instances, it turns out that the order of the variance of 
$LC_n$ is thus linear in $n$ which is the order conjectured 
by Waterman~\cite{WE}, for which 
Steele~\cite{S} had previously obtained a generic linear upper bound.  
However, the most important equiprobable iid (say, binary) case 
remains open.

In the present paper our purpose is different and we study 
optimal alignments at a macroscopic level.  
We present a general methodology showing  
that a property of optimal alignments holds true,  
provided this property typically holds true for short-strings alignments.  
To prove these results, we partition $X$ into pieces of \textit{fixed} 
length $k$ and, as 
$n$ goes to infinity, show that typically, and further with probability 
exponentially close to one,  
in any optimal alignment 
most of these pieces get aligned with 
pieces of $Y$ of similar length $k$.  
In other words, with high probability, 
there are no macroscopic gaps in the optimal alignments.

Let us explain, on a further example, how alignments have 
to stay close to the diagonal.  
Take the two related, English and German, words: $X=mother$ and 
$Y=mutter$. The longest common subsequence is $mter$, hence 
$LC_6=4$ and the common 
subsequence $mter$ corresponds to the following two alignments:
\begin{equation}
\label{firstalignment}
\begin{array}{c|c|c|c|c|c|c|c}
m&o& &t& &h&e&r\\\hline m& &u&t&t& &e&r
\end{array}
\end{equation}
and 
\begin{equation}
\label{secondalignment}
\begin{array}{c|c|c|c|c|c|c|c}
m&o& & &t &h&e&r\\\hline m& &u&t&t& &e&r
\end{array}
\end{equation}
Next, view alignments as subsets of $\mathbb{R}^2$  
as follows:  If the $i$-th letter of $X$ gets aligned with the
$j$-th letter of $Y$, then the set representing the alignment is to
contain $(i,j)$. For example, the alignment \eqref{firstalignment}
can be represented as: $(1,1),(3,3),(5,5),(6,6)$ with the
corresponding plot
\begin{equation}
\label{2Dal1}
\begin{array}{c|c|c|c|c|c|c}
r& & & & & &\bullet\\\hline e& & & & &\bullet& \\\hline t& & & & & & \\\hline t&
& &\bullet& & & \\\hline u& & & & & & \\\hline m&\bullet& & & & & \\\hline
 &m&o&t&h&e&r
\end{array}
\end{equation}
while the alignment \eqref{secondalignment}
can be represented as: $(1,1),(3,4),(5,5),(6,6)$ with the
corresponding plot
\begin{equation}
\label{2Dal2}
\begin{array}{c|c|c|c|c|c|c}
r& & & & & &\bullet\\\hline e& & & & &\bullet& \\\hline t& & &\bullet & & & \\\hline t&
& & & & & \\\hline u& & & & & & \\\hline m&\bullet& & & & & \\\hline
 &m&o&t&h&e&r
\end{array}
\end{equation}
Above, the symbol $\bullet$ indicates the aligned letter-pairs, and 
these points are said to represent the optimal alignment.  
An optimal alignment can then be viewed as the graph 
of a function $f$ defined as follows:  Let $f(0) = 0$ and for $i,j= 1, \dots, n$, 
let $f(i) = j$ if the $ith$ letter 
of the first sequence is aligned with the $jth$ letter of the second sequence, while  
between these values and till $(n,n)$  let 
$f$ be defined via a linear interpolation.  Then as shown in Section 4, 
with its notation, 
the function $f=f_{p_1,p_2}:[0, n]\rightarrow [0,n]$ stays between the 
two lines with 
respective slope $p_1$ and $p_2$.  Moreover, under a strict concavity 
assumption on the limiting shape $\tilde\gamma$ (see the next section), $f_{p_1,p_2}$ becomes 
uniformly close to the identity as $n$ tends to infinity.  More precisely, 
let $g_n:[0, 1]\rightarrow [0,1]$ be defined via $g_n(x) = f(nx)/n$, then 
$\lim_{n\to +\infty}\sup_{0\le x \le 1}|g_n(x) - x| = 0$.  
In other words, there are no macroscopic gaps in any optimal alignment  
and any such alignment must remain close to 
the main diagonal.  

\textit{This closeness to the diagonal property has proved crucial in obtaining the 
first result on the limiting law of $LC_n$, under a lower bound 
on the order of the variance, see \cite{HI}.  Broadly speaking, 
when not close to the diagonal, many terms contribute to our CLT estimation but the 
corresponding set of random alignments has exponentially small probability; while 
when alignments are close to the diagonal, the estimation comes from only a few terms. 
The balance between these 
two cases leads to a central limit theorem.  This CLT contrasts 
with the case of two independent 
uniform random permutations of $\{1, 2, \cdots, n\}$, where the limiting 
distribution of the length of the longest common subsequences is the 
Tracy-Widom distribution, see \cite{HI}.}     
The others, more pathological, instances we are aware of and 
where a CLT holds true in a corresponding LPP problem is for the length of the 
longest increasing subsequences 
of a single random word, potentially Markovian, where only one letter is attained 
with maximal probability, see \cite{ITW1}, \cite{ITW2}, \cite{HL1}, \cite{HL2}. 
(However, for two or more random words the limiting law is no longer Gaussian \cite{BH}.)  
This one word result contrasts, once more, with the single word 
permutation result of Baik, Deift and Johnansson \cite{BDJ} where 
the limiting law is the Tracy-Widom one.

Note that the LCS problem can be represented as a directed 
last passage percolation (LPP) problem with dependent weights.  
Indeed, let the set of vertices  be 
 \[V:=\{0,1,2,\ldots,n\}\times \{0,1,2,\dots,n\},\] 
 and let the set of oriented edges ${\cal E}\subset V\times V$ 
contain horizontal, vertical and diagonal edges.  
The horizontal edges are oriented to the right, 
while the vertical edges are oriented upwards, both having unit length. 
The diagonal edges point up-right at a $45$-degree angle and have length $\sqrt{2}$.  
Hence, 
\[{\cal E} :=\left\{ (v,v+e_1),(v,v+e_2),(v,v+e_3): v\in V\right\},\]
where $e_1:=(1,0)$, $e_2:=(0,1)$ and $e_3:=(1,1)$.  
With the horizontal and vertical edges, we associate 
a weight of $0$.  With the diagonal edge from $(i,j)$ to $(i+1,j+1)$ we associate 
the weight $1$ if $X_{i+1}=Y_{j+1}$ and $0$ (or $-\infty$) otherwise.  
In this manner, we obtain that $LC_n$  
is equal to the total weight of the heaviest paths going from $(0,0)$ to $(n,n)$.  
(Another directed LPP representation can be obtained 
via $LC_n = \max_{\pi \in SI}\sum_{(i,j)\in \pi}{\bf 1}_{\{X_i = Y_j\}}$, 
where $SI$ refers to the set of all paths with \textit{strictly} increasing steps, i.e., paths 
with \textit{both} coordinates strictly increasing from a step to another, from $(0,0)$ to 
the East, $x=n$, or North, $y=n$, boundary.  A third representation would be as above 
but where now the paths 
going from $(0,0)$ to $(n,n)$ have either strictly increasing steps 
or North or East unit steps.  
Again to the strictly 
increasing steps the associated weight is ${\bf 1}_{\{X_i = Y_j\}}$ 
while to the North as well as to the East unit steps is associated a weight value of $0$.  
As a final representation one could still proceed with strictly increasing paths but 
with the requirement that one ends the paths with a $1$.)

Note that the weights in our percolation representations are not ``truly 2-dimensional'' 
and, in our opinion, this is the 
reason for the order of magnitude of the mean, variance as well 
as the limiting law in the LCS problem to be different from 
other first/last passage-related models.  To return to our specific results, 
they are the first studying the transversal fluctuations of the maximal paths of LCSs.  
Such questions have been of much interest in other percolation models.  
Let us only mention that Johansson~\cite{J} showed that for a Poisson points model  
in the plane, typical deviations of a maximal path from the diagonal is of order $n^{2/3}$,  
and $n^{2/3}$ is also the order of the transversal fluctuations 
in the directed polymer model studied in Sepp\"al\"ainen~\cite{Se}.    
(We refer to the respective bibliography of \cite{J} and \cite{Se} for a much more complete 
and detailed picture on transversal fluctuations.  We finally also note that in view of 
the results of \cite{HI}, the transversal fluctuations of the LCSs of two independent  
random permutations of $\{1, 2, \dots, n\}$, one uniform and one arbitrary, are exactly 
the same as those of the longest increasing subsequences of a single 
uniform random permutation 
of $\{1, 2, \dots, n\}$.)

To finish this introductory section, let us briefly describe the rest of the paper.  
The next section presents some preliminary results, examples, and states the main result of 
the paper which is, in turn, proved in Section 3.  Section 4 settles the closeness 
to the diagonal result and shows that each maximal path corresponding to the LCSs stays 
close to the diagonal.  
The last section explores the generic nature of short-strings alignments and some 
of its computational consequences.

\section{Preliminaries}  
\label{sectprelim}

Throughout, let $n=km$   
and let the integers 
\begin{equation}
\label{increasing}
r_0=0<r_1<r_2<r_3<\cdots <r_{m-1}<r_m=n, 
\end{equation}
be such that  
\begin{equation}
\label{sum}
LC_n=
\sum_{i=1}^{m}|LCS(X_{k(i-1)+1}X_{k(i-1)+2}\cdots X_{ki};Y_{r_{i-1}+1}Y_{r_{i-1}+2}
\cdots Y_{r_{i}})|, 
\end{equation}
where $|LCS( \cdot ; \cdot )|$ is the length of the corresponding longest common substrings.  
In words, \eqref{sum} asserts that there exists an optimal alignment 
aligning $[k(i-1)+1,ki]_{\bbn}$ with $[r_{i-1}+1,r_{i}]_{\bbn}$, 
for all $i=1,2,\dots,m$.

The first goal of the present paper is to show that 
for $k$ large but fixed, and $n$ large enough,  
any such generic optimal alignment is such that the vast majority 
of the intervals $[r_{i-1}+1,r_{i}]_{\bbn}$ have length close to $k$. 
Building on this, a second goal is to show 
(see Section~\ref{sectionproperty}) that if a property 
$\mathcal{P}$ holds, with high probability, for string-pairs of (short) length order $k$, then 
typically a large proportion 
of aligned string-pairs satisfy the property $\mathcal{P}$.

Let us deal with our first goal and show that with high probability the 
optimal alignments satisfying \eqref{sum}, are such that most of 
their lengths $r_{i}-r_{i-1}$ 
are close to $k$.   Of course, we need to quantify what is meant 
by ``close to $k$".  To 
do so, we first provide a definition.  For $p>0$, let 
\begin{equation}\label{defgam}
\tilde\gamma(p):=\lim_{n\rightarrow\infty}
\frac{\bbe |LCS(X_1X_2\cdots X_{n};Y_1Y_2\cdots Y_{np})|}{n(1+p)/2}, 
\end{equation} 
where, when not integers, the indices $np$ 
are understood to be rounded-up to the nearest 
positive integers.
This function $\tilde\gamma$ is just a re-parametrization 
of ``the usual function"  
\[\gamma(q)=\lim_{n\rightarrow\infty}
\frac{\bbe |LCS(X_1X_2\cdots X_{n-nq};Y_1Y_2\cdots Y_{n+nq})|}{n},
\]
$q\in(-1,1)$, i.e., 

\begin{equation}
\label{gammagamma*}
\tilde\gamma(p)=\gamma\left(q(p)\right), 
\end{equation}
with $q(p) = (p-1)/(p+1)$

A superadditivity argument, as in Chv\'atal and Sankoff \cite{CS},  
shows that the above limits do exist (and depend,
for example, on the size of the alphabet but this is of no importance 
for our purposes).  For $X$ and $Y$ identically distributed, 
the function $\gamma$ is symmetric about 
the origin, while a further superadditivity argument
shows that it is concave there and so it reaches its maximum 
at $q=0$ (see \cite{AHM1} for details).   It is not known whether or not 
the function $\gamma$ 
is strictly concave around $q=0$.   From 
simulations this appears to be the case but, at present, a proof is elusive. 
(Again, the LCS problem is a last passage percolation problem, and for 
first/last passage percolation proving that 
the shape of the wet region is strictly concave seems difficult and in many cases has 
not been done.)  Since $q(p)=(p-1)/(p+1)=1-2/(p+1)$  
is strictly increasing in $p$, with $q(1)=0$, 
if $\gamma$ were strictly concave around $q=0$, then it would reach 
a strict maximum there. 
Thus, $\tilde\gamma$ would also reach a strict maximum 
at $p=1$ but, without the strict concavity of $\gamma$,   
$p=1$ might not be the unique point of maximal value.

However, strict concavity is not needed, concavity is enough, 
for our results to hold.  Indeed,  
(see Lemma~\ref{expdiff}) $\tilde\gamma$ is non-decreasing 
on $[0,1]$ and non-increasing on $[1,\infty)$ and 
thus ad hoc methods will show that 
$\tilde\gamma(p)$ is strictly smaller than 
$\tilde\gamma(1) = \gamma^*$, as soon as $p$ is farther away from $1$ than a small given 
quantity. (If the unproven strict concavity is valid 
then, below, $p_1$ and $p_2$ could be chosen as close to $1$ as one wishes to.)

Before coming to the proof of the main theorem, it is 
thus important to show the existence of, and provide 
estimates on, the aforementioned $p_1$ and $p_2$ chosen as close to $1$ as possible.  
Let us first explore this question in the binary equiprobable case.  
Let $\gamma_\ell$ be any strict lower bound on  
$\gamma^*=\lim_{n\rightarrow\infty} \bbe|LCS(X_1\cdots X_n;Y_1\cdots Y_n)|/n$.  
Then, if $x>0$ is such that  
\begin{eqnarray}
\label{entropy2}
&&\!\!\!\!\!\!\!(1-x)\!\left(\!H_2\!\left(\!1-\frac{\gamma_\ell(2-x)}{2(1-x)}\!\right)-1\!\right)
+H_2\!\left(\!x+(1-x)\left(\!1-\frac{\gamma_\ell(2-x)}{2(1-x)}\!\right)\!\right) \nonumber\\
&&= (1-x)\!\left(\!H_2\!\left(\!1-\frac{\gamma_\ell}{2}\!\left(1+\frac{1}{1-x}\right)\!\right)-1\!\right)
+H_2\!\left(\!1-\frac{\gamma_\ell}{2}\!\left(1+(1-x) \right)\!\right)<0, 
\end{eqnarray}
where $H_2(x) = -x\log_2 x - (1-x)\log_2(1-x)$, $0<x<1$, 
is the binary entropy function,  we claim that 
$p_1:=1-x$ and $p_2:=1/(1-x)$ are such that 
\begin{equation}\label{gammap1p2zero}
\tilde\gamma(p_1)<\tilde\gamma(1)=
\gamma^*,\;\;\;\;\tilde\gamma(p_2)<\tilde\gamma(1) = \gamma^*.     
\end{equation}
Indeed, an easy upper bound on the probability that the length of the LCSs of $X_1\cdots X_n$ 
and $Y_1\cdots Y_{n(1-x)}$  
is larger than $n(1-\eta)(1-x)$, $0<\eta<1$ is found as follows:  
First, take a non-random string $s$ 
of length $n(1-\eta)(1-x)$.  The probability that $s$ is a 
subsequence of $Y_1\cdots Y_{(1-x)n}$ 
does not depend on $s$ (but only on the length of $s$), and so this probability 
would be the same if $s$ would consist only of ones.  
Therefore, the probability that $s$ is a subsequence of 
$Y_1\cdots Y_{n(1-x)}$ is the same as the probability that there are at least 
$n(1-x)(1-\eta)$ ones in $Y_1\cdots Y_{n(1-x)}$, which, in turn, is nothing 
but the probability 
for a binomial random variable, with parameters $(1-x)n$ and $1/2$, 
to be greater or equal to $n(1-x)(1-\eta)$.  
But, via classical exponential inequalities, this last probability is bounded above by 
\begin{equation}
\label{thatha}
2^{n(1-x)(H_2(\eta)-1)}.
\end{equation}
Now, the number of subsequences of $X_1\cdots X_n$ of length  $n(1-\eta)(1-x)$ is given by:  
\[\binom{n}{n(\eta+x-\eta x)}\le 2^{nH_2(\eta +x-x\eta)}.\]
Combining this last inequality with \eqref{thatha} leads to
\begin{equation}
\label{name}
\bbp\left(|LCS(X_1\cdots X_n;Y_1\cdots Y_{n(1-x)})|\geq (1-x)(1-\eta) n\right)
\le 2^{n((1-x)(H_2(\eta)-1)+H_2(\eta+x-\eta x))}.  
\end{equation}
Therefore, from \eqref{name}, as soon as 
\begin{equation}
\label{always}
(1-x)(H_2(\eta)-1)+H_2(\eta+x-\eta x)<0,
\end{equation}
then, the probability that $|LCS(X_1\cdots X_n;Y_1\cdots Y_{n(1-x)})|$ 
is at least $(1-x)(1-\eta)n$ is exponentially
small in $n$.   In other words, the probability that the rescaled 
(by the average length of the two
strings) LCS is at least:
\begin{equation}
\label{thotho}
\frac{2(1-x)(1-\eta)}{2-x}, 
\end{equation}
is exponentially small in $n$.  
Now, to require the rescaled LCS value to be equal to $\gamma_\ell$, 
choose $\gamma_\ell$ equal to the quantity in \eqref{thotho}, i.e., let 
\begin{equation}
\label{heiniheini}
\eta=1-\frac{\gamma_\ell(2-x)}{2(1-x)}.
\end{equation}
So, for $x$ given, \eqref{heiniheini} 
gives the value of $\eta$ corresponding to $\gamma_\ell$. 
We next find values of $x$ for which the probability,  
that the rescaled LCS of $X_1\cdots X_n$ and $Y_1\cdots Y_{n(1-x)}$  
is at least $\gamma_\ell$, is exponentially small in $n$.  
Indeed, in \eqref{always} it is enough to replace $\eta$ by 
the value given in \eqref{heiniheini}.  
This leads to the condition \eqref{entropy2}, and when 
\eqref{entropy2} is satisfied, the probability 
that the LCS of $X_1\cdots X_n$ and $Y_1\cdots Y_{n(1-x)}$ 
has a rescaled value of at least $\gamma_\ell$ is exponentially small in $n$.  
Therefore, in this case, 
the rescaled limit, that is 
\[\tilde\gamma(1-x)
=\lim_{n\rightarrow\infty}\frac{\bbe|LCS(X_1\cdots X_n;Y_1\cdots Y_{n(1-x)})|}{n},\] 
has to be 
at most $\gamma_\ell$.  
Thus, in this case, if $\gamma_\ell$ is a lower bound on $\tilde\gamma(1)$, 
then \[\tilde\gamma(1-x)\le \gamma_\ell < \tilde\gamma(1)=\gamma^*,\]
so that $p_1=1-x$ and $p_2=1/(1-x)$ satisfy \eqref{gammap1p2zero}.

\

Using $\gamma_\ell = 0.7880$, which is a lower bound on $\gamma^*$ 
obtained, in the binary case, by Lueker~\cite{L}, 
it is easily seen that \eqref{entropy2} is negative for $x = 0.28$.  
So $p_1 = 0.72$ and $p_2 = 1.39$ satisfy the needed conditions.  
Using these values as well as the upper bound 
$0.8263$ (see \cite{L}) provide an estimate for 
the fixed length $k$ for the pieces one would divide the strings into 
(see, below, the statement of our first theorem).

\

The entropic method, on obtaining bound on $p_1$ and $p_2$, presented above 
carries over, beyond the binary case, to arbitrary-size 
alphabets.  However, in the non-uniform case such bounds might be far from optimal 
and require the knowledge of the probability associated to each letter.  To further 
address this question, let us present a lemma, with a somehow easier 
approach, to deal with the generic case.    

\

\begin{lemma}\label{generalpestimates} 
Let $0< p_1 = \gamma^*/(2-\gamma^*) < 1 $ and let $1< p_2 = (2-\gamma^*)/\gamma^*$, then   
\begin{equation}\label{gammap1p2two}
\tilde\gamma(p_1)<\tilde\gamma(1)=
\gamma^*,\;\;\;\;\tilde\gamma(p_2)<\tilde\gamma(1) = \gamma^*.     
\end{equation}
Moreover, \eqref{gammap1p2two} continue to hold by taking 
$p_1 = \gamma_\ell/(2-\gamma_\ell)$ and 
$p_2 = (2-\gamma_\ell)/\gamma_\ell$ where $\gamma_\ell$ is any positive lower bound,    
such as  $\sum_{\alpha\in\cala}(\bbp(X_1=\alpha))^2$, on $\gamma^*$.   
\end{lemma}

\begin{proof}
First, note that when one sequence 
has length $0$, then the LCS also has length $0$, and thus 
\[\gamma(-1)=\gamma(1)=0.\]  
Recall next that the function $\gamma$ is concave  
and symmetric about $0$.  Next, consider the LCS of the string 
$X=X_1X_2\cdots X_{2n}$ and the empty string $Y$.  Then, take 
$\epsilon n$ letters away from $X_1\cdots X_{2n}$ and add 
that many letters to the $Y$ string, so as to now have,   
instead of the empty string, the string $Y_1Y_2\cdots Y_{\epsilon n}$.
Then, provided $0<\epsilon<1$ is a small enough constant, it follows 
that with very high (in $n$) probability, 
the string $Y_1Y_2\cdots Y_{\epsilon n}$ is a substring of 
$X_1\cdots X_{2n-\epsilon n}$.   Indeed, let $T_1=\inf\{k\ge 0: X_1 = Y_k\}$, 
$T_2=\inf\{k\ge T_1: X_2 = Y_k\}$, $T_3=\inf\{k\ge T_2: X_3 = Y_k\}, \dots$  
Then, conditionally on $X=x$, $T_1, T_2-T_1, T_3-T_2, \dots$ are independent 
geometric random variables with individual parameter depending on the 
probabilities associated to the letters.  But, 
\[|LCS(X_1\cdots X_{2n-\epsilon n};Y_1\cdots Y_{\epsilon n})|=\epsilon n,\] 
if and only if $T_{\epsilon n} \le 2n -\epsilon n$.  Thus, since the geometric 
property holds for any $x\in{\cala}^n$,  
\begin{equation}
\label{gatech}
\bbp(|LCS(X_1\cdots X_{2n-\epsilon n};Y_1\cdots Y_{\epsilon n}|)=\epsilon n)\ge 1-\exp(-cn), 
\end{equation}
where $c>0$ is a constant depending neither on $n$ nor on $\epsilon>0$ (provided 
$\epsilon>0$ is small enough) but depending on the minimal parameter of the geometric 
random variables (hence on the probabilities associated to the letters).  
Therefore, from the above, at $q=-1$, the slope is $1$, i.e., 
$\gamma^\prime((-1)^+)=1$ and similarly $\gamma^\prime(1^-)=-1$.  
By concavity and symmetry, for any $q_1$ with 
\begin{equation}
\label{conditionq_1}
q_1<-(1-\gamma(0)), 
\end{equation} 
it therefore follows that 
\[\gamma(q_1)<\gamma(0),\]
and similarly for any $q_2$ with 
\begin{equation}
\label{conditionq_2}
q_2>1-\gamma(0), 
\end{equation}
it follows that 
\[\gamma(q_2)<\gamma(0).\]
The bounds obtained above rely on the value of $\gamma(0)= \gamma^*$ 
which is unknown, but for which upper and lower bounds 
(which depend on the distributions of the letters) 
do exist (and are rather accurate for uniform distributions).  
In our case, a lower bound on $\gamma^*$ is what is needed.  
The most trivial 
lower bound is obtained when aligning the two strings, without gaps, 
and just counting the number of correctly aligned letter pairs, more precisely, 
$LC_n\ge \sum_{i=1}^n{\bf 1}_{\{X_i=Y_i\}}$.   
Hence, by the iid and independence assumptions, 
\begin{equation}
\label{lower}\bbe LC_n\ge\sum_{i=1}^n\bbp(X_i=Y_i)
=n\bbp(X_1=Y_1)=n\sum_{\alpha\in \cala}\bbp((X_1=\alpha))^2.
\end{equation}
Now, converting to 
$\tilde{\gamma}$, and recalling that $p(q)=(1+q)/(1-q)$, 
\eqref{conditionq_1} becomes 
\begin{equation}
\label{boundp1bis}
p_1<\frac{\gamma(0)}{2-\gamma(0)}, 
\end{equation}
and in such a case, 
$\tilde{\gamma}(p_1)<\tilde{\gamma}(1)$.  
Similarly, \eqref{conditionq_2}
becomes 
\begin{equation}
\label{bound2bis}
p_2>\frac{2-\gamma(0)}{\gamma(0)}, 
\end{equation}
and then
$\tilde{\gamma}(p_2)<\tilde{\gamma}(1)$.  
Finally, in both \eqref{boundp1bis} and \eqref{bound2bis}, 
one can replace $\gamma(0)$
by any of its positive lower bound, for example the lower 
bound resulting from \eqref{lower}. \end{proof}

So let us now assume that
$0 < p_1 < 1 < p_2$ are such that 
\begin{equation}
\tilde\gamma(p_1)<\tilde\gamma(1)=
\gamma^*,\;\;\;\;\tilde\gamma(p_2)<\tilde\gamma(1) = \gamma^*.   
\end{equation}
The first result we next set to state, asserts that for $k$ 
fixed and $n$ large enough 
then typically in any optimal alignment \eqref{sum}  
most of the intervals $[r_{i-1}+1, r_i]_{\bbn}$ (for $i=1,2, \dots, m$) 
have their length, $r_i-r_{i-1}$, between $kp_1$ and $kp_2$. By most, 
it is meant that by taking $k$ large, \textit{but fixed}, and $n$ large enough, 
the proportion of such intervals gets as close to $1$ as one wishes to.   

To do so, let us introduce some notation.  
Let $\epsilon>0$, $p_1> 0$ and $p_2>0$ be constants.
Let $A^n_{\epsilon,p_1,p_2}$ 
be the (random) set of optimal alignments 
of $X_1\cdots X_n$ and $Y_1\cdots Y_n$ satisfying \eqref{sum}, 
for which a proportion of at least $1-\epsilon$ 
of the intervals $[r_{i-1}+1,r_{i}]_{\bbn}$, $i=1,2,\dots,m$, 
have their length between $kp_1$ and $kp_2$.  
More precisely, $A^n_{\epsilon,p_1,p_2}$ is the event that for 
all integer vectors $(r_0,r_1,\dots,r_m)$
satisfying \eqref{increasing} and for which \eqref{sum} holds, 
\begin{equation}
\label{condition2} {\rm Card}
\left(\;\{i\in 1,2,\dots, m: kp_1\le r_i-r_{i-1}\le kp_2\}\;\right) 
\geq (1-\epsilon)m.
\end{equation}

As stated next, $A^n_{\epsilon,p_1,p_2}$ holds with high probability.    

\begin{theorem}
\label{maintheorem}
Let $\epsilon>0$. Let $0< p_1<1<p_2$ be such that 
$\tilde\gamma(p_1)< \tilde\gamma(1)=\gamma^*$ and $\tilde\gamma(p_2)
<\tilde\gamma(1)=\gamma^*$,  
and let 
$\delta\in (0,\min(\gamma^*-\tilde\gamma(p_1),\gamma^*-\tilde\gamma(p_2)))$.  
Fix the integer $k$ to be such that $(1 + \ln k)/k\le \delta^2\epsilon^2/16$, then 
\begin{equation}\label{mainineq}
\bbp\left(A^n_{\epsilon,p_1,p_2}\right)\geq 1-
\exp\left(-n \left(-\frac{1+\ln k}{k}+\frac{\delta^2\epsilon^2}{16}\right)\right),
\end{equation}
for all $n = n(\epsilon, \delta)$ large enough. 
\end{theorem}
 
In words, and broadly, Theorem~\ref{maintheorem} asserts that for any $\epsilon > 0$, 
there exists 
$k$ large enough, but fixed, such that if $X$ is divided into segments of length $k$ then,  
typically (at least a fraction $1-\epsilon$ of segments), and with high probability, 
the LCSs match these segments to segments of similar length in $Y$.

\section{Proof of the main theorem}
\label{sectmainth}

The proof of Theorem~\ref{maintheorem} requires the introduction of 
a few more definitions.  So far we have looked at the integer intervals 
which are mapped by optimal alignments 
to the integer intervals $[k(i-1)+1, ki]_{\bbn}$.  The opposite stand is now taken:  
given (non-random) integers 
$r_0=0<r_1<r_2<\cdots<r_m=n$, we request that the alignment 
aligns $[k(i-1)+1,ki]_{\bbn}$ with $[r_{i-1}+1,r_i]_{\bbn}$, for every $i=1,2,\dots,m$. 
In general, such an alignment is not optimal and, the best 
score an alignment can attain under the above constraint is:
\begin{align*}
L_n(\vec{r})&:=L_n(r_0,r_1,\dots,r_m)\\
&:
=\sum_{i=1}^{m}|LCS(X_{k(i-1)+1}X_{k(i-1)+2}\cdots X_{ki};Y_{r_{i-1}+1}Y_{r_{i-1}+2}\cdots Y_{r_{i}})|.
\end{align*}
Therefore, $L_n(\vec{r})$ represents the maximum 
number of aligned identical letter pairs under the constraint 
that $X_{(i-1)k+1}X_{(i-1)k+2}\cdots X_{ik}$ gets aligned with 
$Y_{r_{i-1}+1}Y_{r_{i-1}+2}\cdots Y_{r_i}$, for all $i=1,2,\dots, m$.  
Note moreover that for a non-random $(r_0,r_1,\dots,r_m)$, the 
partial scores 
\[|LCS(X_{(i-1)k+1}X_{(i-1)k+2}\cdots X_{ik};
Y_{r_{i-1}+1}Y_{r_{i-1}+2}\cdots Y_{r_i})|,\]
are independent of each other and, in this context, concentration inequalities will 
prove useful when dealing with $L_n(\vec{r})$.  Next, let $\mathcal{R}_{\epsilon,p_1,p_2}$, 
be the (non-random) set of all integer vectors 
$(r_0,r_1,\dots,r_m)$ satisfying \eqref{increasing} 
and \eqref{condition2}, while $\overline{\mathcal{R}}_{\epsilon,p_1,p_2}$, 
denotes the (non-random) set of all integer vectors
$\vec{r}=(r_0,r_1,\dots,r_m)$ satisfying \eqref{increasing} but not \eqref{condition2}.

\

Let us begin with a lemma.

\begin{lemma}\label{expdiff}
Let $\epsilon > 0$.  Let $0<p_1<1<p_2$ be such that $\tilde\gamma(p_1)<\gamma^*$ and 
$\tilde\gamma(p_2)<\gamma^*$, and  let $\delta>0$ be such that 
$\delta<\min(\gamma^*-\tilde\gamma(p_1),\gamma^*-\tilde\gamma(p_2))$.  Let 
$\vec{r}=(r_0,\dots,r_m)\in\overline{\mathcal{R}}_{\epsilon,p_1,p_2}$,   
then  
\begin{equation}
\label{gammafrac12}
\bbe\left(L_n(\vec{r})-LC_n\right)\le -\frac{\delta \epsilon n}{2},   
\end{equation}
for all $n = n(\epsilon, \delta)$ large enough.
\end{lemma}

\begin{proof} 
Let $p>0$, and let  
\[\tilde\gamma(p):=
\lim_{n\rightarrow\infty}\frac{\bbe|LCS(X_1X_2\cdots X_n;Y_1Y_2\cdots Y_{np})|}{n(1+p)/2}, \]
where, by superadditivity, this limit exists with, moreover,  
\begin{equation}
\label{notatall}
\frac{2\bbe|LCS(X_1X_2\cdots X_n;Y_1Y_2\cdots Y_{np})|}{n(1+p)}\le \tilde\gamma(p),
\end{equation}
for any $n\geq 1$.  Now, $\gamma: q\in (-1,1)\rightarrow\gamma(q)\in (0,\infty)$, defined via
\[\gamma(q):=\lim_{n\rightarrow\infty}
\frac{\bbe|LCS(X_1X_2\cdots X_{n-nq};Y_1Y_2\cdots Y_{n+nq})|}{n},\]
is symmetric about $q=0$ and, as already mentioned, is also concave (see \cite{AHM1}).  
Hence,  
\[\tilde\gamma(p)=\gamma\left(\frac{p-1}{p+1}\right),\]  
is non-decreasing up to $p=1$ and non-increasing afterwards. 
Thus, choosing the interval $[p_1,p_2]$ to contain $p=1$, it follows   
that for $p \notin [p_1,p_2]$, 
\begin{equation}
\label{myrta}
\tilde\gamma(p)\le \max(\tilde\gamma(p_1),\tilde\gamma(p_2)).
\end{equation}

Hence, for any $p \notin [p_1,p_2]$, combining \eqref{notatall} 
and \eqref{myrta} leads to:  
\begin{equation}
\label{samuel0}
\frac{2\bbe|LCS(X_1X_2\cdots X_k;Y_1Y_2\cdots Y_{kp})|}{k(1+p)}\le 
\max(\tilde\gamma(p_1),\tilde\gamma(p_2)),    
\end{equation}
and therefore, 
\begin{equation}
\label{samuel20}
\frac{2\bbe|LCS(X_1X_2\cdots X_k;Y_1Y_2\cdots Y_{kp})|}{k(1+p)}\le 
\tilde\gamma(1)-\delta^* = \gamma^* - \delta^*, 
\end{equation}
where $\delta^*:=\min(\gamma^*-\tilde\gamma(p_1),\gamma^*-\tilde\gamma(p_2))$.

Since the sequences $(X_i)_{i\ge 1}$ and $(Y_i)_{i\ge 1}$ are stationary,  
and assuming that $r_i-r_{i-1}=kp$, 
the left-hand side of \eqref{samuel20} becomes 
\[\frac{2\bbe|LCS(X_{(i-1)k+1}X_{(i-1)k+2}\cdots X_{ik};
Y_{r_{i-1}+1}Y_{r_{i-1}+2}\cdots Y_{r_i})|}{k+r_i-r_{i-1}}.\]
Thus, from \eqref{samuel20}, when $(r_i-r_{i-1})/k \notin [p_1,p_2]$, then 
\begin{equation}\label{encoreun}
\gamma^*-\frac{2\bbe|LCS(X_{(i-1)k+1}X_{(i-1)k+2}\cdots X_{ik};
Y_{r_{i-1}+1}Y_{r_{i-1}+2}\cdots Y_{r_i})|}{k+r_i-r_{i-1}}\geq\delta^*.
\end{equation}

Hence, from \eqref{encoreun}, 
\begin{align*}
\gamma^*\left(\frac{k+r_i-r_{i-1}}{2}\right) &- \bbe|LCS(X_{(i-1)k+1}X_{(i-1)k+2}\cdots X_{ik};
Y_{r_{i-1}+1}Y_{r_{i-1}+2}\cdots Y_{r_i})|  \\ 
&\qquad \qquad \ge \delta^*\left(\frac{k+r_i-r_{i-1}}{2}\right) \ge 
\delta^* \frac{k}{2}.
\end{align*}
Letting ${\cal M} := \{i:[k(i-1)+1, ki]$
gets  matched  with  strings  of  length not  in  $[kp_1, kp_2]\}$, 
we then have  
\begin{align}\label{encoredeux}
\sum_{i\in {\cal M}}\Big(\gamma^*\Big(\frac{k+r_i-r_{i-1}}{2}\Big) 
& - \bbe|LCS(X_{(i-1)k+1}X_{(i-1)k+2}\cdots X_{ik};
Y_{r_{i-1}+1}Y_{r_{i-1}+2}\cdots Y_{r_i})|  \Big) \nonumber \\
&\ge \sum_{i\in {\cal M}}\delta^*\frac{k}{2} 
\ge \delta^*\frac{k}{2}\epsilon m = \frac{n\delta^*\epsilon}{2}.
\end{align}
On the other hand, 
\begin{align}\label{encoretrois}
&\sum_{i\in {\cal M}}\Big(\gamma^*\Big(\frac{k+r_i-r_{i-1}}{2}\Big) 
- \bbe|LCS(X_{(i-1)k+1}X_{(i-1)k+2}\cdots X_{ik};Y_{r_{i-1}+1}Y_{r_{i-1}+2}
\cdots Y_{r_i})|\Big) \nonumber\\ 
&\le \sum_{i=1}^m \Big(\gamma^*\Big(\frac{k+r_i-r_{i-1}}{2}\Big) 
- \bbe|LCS(X_{(i-1)k+1}X_{(i-1)k+2}\cdots X_{ik};Y_{r_{i-1}+1}Y_{r_{i-1}+2}\cdots Y_{r_i})|\Big)
\nonumber \\
&= \gamma^*n - \bbe L_n(\vec{r}).
\end{align}

\noindent
Therefore, combining \eqref{encoredeux} and \eqref{encoretrois} leads to  
\begin{equation}
\label{star1}
\gamma^* n-\bbe L_n(\vec{r}) \geq \frac{n\delta^*\epsilon}2,
\end{equation}
as soon as $(r_0,r_1,\dots,r_m)\in \overline{\mathcal{R}}_{\epsilon,p_1,p_2}$.  
Now $\lim_{n\to +\infty}\bbe LC_n/n = \gamma^*$, while by hypothesis $\delta^*-\delta>0$, 
\begin{equation}
\label{star2}
0 \le \gamma^*-\frac{\bbe LC_n}{n}\le \frac{(\delta^*-\delta)\epsilon}{2}.  
\end{equation}
for all $n$ large enough.  
(At this last stage, a more quantitative bound, depending on $n$, could also 
be obtained using \eqref{alex}.) 
Combining \eqref{star1} and \eqref{star2} yields that for any 
$\vec{r}\in\overline{\mathcal{R}}_{\epsilon,p_1,p_2}$,  
\[\bbe\left(LC_n-L_n(\vec{r})\right)\geq \frac{n\delta\epsilon}{2},
\]
for all $n$ large enough.
The proof of the lemma is now complete.\end{proof}

\begin{proof}[Proof of Theorem~\ref{maintheorem}]
Clearly, 
\begin{equation}\label{inebin}
{\rm Card}\left(\,\overline{\mathcal{R}}_{\epsilon,p_1,p_2}\right)
\le \binom {n}{m} \le \frac{n^m}{m!}\le \left(\frac{en}{m}\right)^m = (ek)^m, 
\end{equation}
by a well known and simple bound on the binomial coefficients.

Now, let $\delta<\delta^*:=\min(\gamma^*-\tilde\gamma(p_1),\gamma^*-\tilde\gamma(p_2))$.  
By definition $LC_n \ge L_n(\vec{r})$, and so for $\vec{r}$ to define an 
optimal alignment requires:
 \begin{equation}
\label{LnLCn}
L_n(\vec{r})\ge LC_n.
\end{equation} 
Hence, for the event 
$A^n_{\epsilon,p_1,p_2}$ not to hold (see \eqref{condition2}), 
there needs to exist at least one 
$\vec{r} \in \overline{\mathcal{R}}_{\epsilon,p_1,p_2}$ 
for which \eqref{LnLCn} is satisfied. Thus, 
\[(A^{n}_{\epsilon,p_1,p_2})^c=
\bigcup_{\vec{r}\in \overline{\mathcal{R}}_{\epsilon,p_1,p_2} }\{L_n(\vec{r})-LC_n\geq 0\},\]
and 
\begin{equation}
\label{samuel}
\bbp((A^{n}_{\epsilon,p_1,p_2})^c)\le
\sum_{\vec{r}\in \overline{\mathcal{R}}_{\epsilon,p_1,p_2} }\bbp(L_n(\vec{r})-LC_n\geq 0).
\end{equation}
When $\vec{r}\in \overline{\mathcal{R}}_{\epsilon,p_1,p_2}$, it follows from 
Lemma~\ref{expdiff} that:  
\[\bbe\left(L_n(\vec{r})-LC_n\right)\le -\frac{\delta \epsilon n}{2},\]
and so 
\begin{equation}
\label{together}
\bbp(L_n(\vec{r})-LC_n\geq 0)\le
\bbp\left(L_n(\vec{r})-LC_n-\bbe\left(L_n(\vec{r})-LC_n\right)\ge 
\frac{\delta\epsilon n}2 \right), 
\end{equation}
for all $n$ large enough.  Now, the difference $L_n(\vec{r})-LC_n$ changes by at most two, 
when any one of the iid entries 
$X_1, X_2, \dots, X_n, Y_1, Y_2,\dots, Y_n$ is changed.  Therefore,  
Hoeffding's martingale inequality, applied to the right-hand side 
of \eqref{together}, gives
\[
\bbp\left(L_n(\vec{r})-LC_n\geq 0\right)\le
\bbp\left(L_n(\vec{r})-LC_n\;-\;\bbe\left(L_n(\vec{r})-LC_n\right) 
\geq \frac{\delta\epsilon n}{2}\right)
\le \exp\left(-\frac{\delta^2 \epsilon^2}{16}n\right). 
\]
(Recall that Hoeffding's martingale inequality asserts that if $f$ is a 
function of $j$ variables, such that changing any single of its entries  
changes $f$ by at most $\Delta$, and if $Z_1, Z_2,\dots, Z_j$ are 
independent random variables, then 
\[\bbp\left(f(Z_1,Z_2,\dots,Z_j)-
\bbe f(Z_1,Z_2,\dots,Z_j)\ge z \right)
\le \exp\left(-\frac{2z^2}{j \Delta^2}\right),\]
provided the expectation exists.)  
Combining this last inequality with \eqref{samuel}, 
one obtains: 
\begin{equation}
\label{samuel2}
\bbp((A^{n}_{\epsilon,p_1,p_2})^c)\le 
\mbox{Card}( \overline{\mathcal{R}}_{\epsilon,p_1,p_2}) 
\exp\left(-\frac{\delta^2 \epsilon^2}{16}n\right). 
\end{equation}
But, from \eqref{inebin}, 
\begin{equation}
\label{samuel3}
\bbp((A^{n}_{\epsilon,p_1,p_2})^c)\le
(ek)^m
 \exp\left(-\frac{\delta^2 \epsilon^2}{16}n\right)=
\exp\left(-n \left(-\frac{1+\ln k}{k}+\frac{\delta^2\epsilon^2}{16}\right)\right). 
\end{equation}
Therefore, the proof of Theorem~\ref{maintheorem} is complete.\end{proof}

\section{Closeness to the diagonal}

Let us begin with a definition.  Let $D^n_{\epsilon,p_1,p_2}$ be the 
event that all the points representing any optimal alignment
of $X_1X_2\cdots X_n$ with $Y_1Y_2\cdots Y_n$
are above the line $y=p_1x-p_1n\epsilon-p_1k$, and below
the line $y=p_2x+p_2n\epsilon+p_2k$.
\begin{theorem}\label{diagol}
Let $\epsilon >0$.  
Let $0< p_1<1<p_2$ be such that $\tilde\gamma(p_1)< \gamma^*$ and
$\tilde\gamma(p_2)<\gamma^*$, and let 
$0< \delta < \min\left(\gamma^*-\tilde\gamma(p_1), \gamma^*-\tilde\gamma(p_2)\right)$.  
Fix the integer $k$ to be such that $(1+\ln k)/k\le \delta^2\epsilon^2/16$, then
\begin{equation}\label{double}
\bbp(D^{n}_{\epsilon,p_1,p_2})\geq 
1-2\exp\left(-n\left(-\frac{1+\ln k}{k}+\frac{\delta^2\epsilon^2}{16}\right)\right),
\end{equation}
for all $n=n(\epsilon,\delta)$ large enough.
\end{theorem}

\begin{proof}
Let $D^n_a$ be the event that any optimal alignment 
of  $X_1X_2\cdots X_n$ with $Y_1Y_2\cdots Y_n$ 
is above the line $y_1:=y_1(x)=p_1x-p_1n\epsilon-p_1k$; and let $D^n_b$ be 
the event that any optimal alignment 
of $X_1X_2\cdots X_n$ with $Y_1Y_2\cdots Y_n$
is below
the line $y_2:=y_2(x)=p_2x+p_2n\epsilon+p_2k$.  Clearly,  
$D^n_a\cap D^n_b = D^n_{\epsilon,p_ 1,p_2}$, hence   
\begin{equation}
\label{bbp}
\bbp((D^{n}_{\epsilon,p_1,p_2})^c)\le 
\mathbb{P}((D^{n}_a)^c) + \mathbb{P}((D^{n}_b)^c),  
\end{equation}
and the result will be a consequence of the following two inclusions:  
\begin{equation}
\label{AnDna}
A^n_{\epsilon,p_1,p_2}\subset D^n_a, \quad A^n_{\epsilon,p_1,p_2}\subset D^n_b,   
\end{equation}
where $A^n_{\epsilon,p_1,p_2}$ is as in Theorem~\ref{maintheorem}.  
Let us prove the first inclusion in \eqref{AnDna}.  
To start, assume that $x$ is an integer multiple of $k$, i.e., 
let $x=uk$, $u\in \bbn$.  Next, and at 
first, let us consider the case where $x\le n\epsilon$, i.e., that 
$p_1x-p_1n\epsilon \le 0$.  Any alignment (and, in particular, any optimal 
alignment) we consider, aligns any $x\in[0,n]_{\bbn}$ with $[0,n]_{\bbn}$.  
Hence, for every $x\le n\epsilon$, the condition is always verified, 
that is any optimal alignment aligns $x$ with a $y$ 
which is at least equal to $p_1x-p_1n\epsilon$.  
Let us now consider the case where $x\ge n\epsilon$.  
When the event $A^n_{\epsilon,p_1,p_2}$ holds, 
then any optimal alignment aligns all but a proportion 
$\epsilon$ of the interval $[(i-1)k+1,ik]_{\bbn}$, $i\in\{1,\dots,m\}$  
to integer intervals of length greater or equal to $kp_1$.  
The maximum number of integer intervals 
which could be matched with integer intervals of length less than 
$kp_1$ is thus $\epsilon m$.  
In the interval $[0,x]_{\bbn}$ there are $u$ intervals 
from the partition $[(i-1)k+1,ik]_{\bbn}$, $i\in\{1,\dots,m\}$.  
Therefore, at least $u-\epsilon m$ of these intervals are matched 
to intervals of length no less than $kp_1$, implying 
that, when the event $A^n_{\epsilon,p_1,p_2}$ holds, 
$x$ gets matched by the optimal alignment 
to a value no less than $(u-\epsilon m)kp_1 = p_1x-p_1\epsilon n$,  
since $x=uk$ and $n=mk$.  This finishes the case where $x$ is an integer 
multiple of $k$.   If $x$ is not an integer multiple of 
$k$, let $x_1$ denote the largest
integer multiple of $k$ which is smaller than $x$.  
By definition,  
\begin{equation}
\label{xx1k}
x-x_1< k.
\end{equation} 
But, the two-dimensional alignment curve cannot 
go down, hence $x$ gets aligned with a point which cannot be  
below the point where $x_1$ gets aligned to.  But, since $x_1$ is an integer multiple of $k$,  
it gets aligned to a point which is greater or equal to $p_1x_1-p_1\epsilon n$.  
Using \eqref{xx1k}, it follows that 
\[p_1x_1-p_1\epsilon n\ge p_1x-p_1\epsilon n-p_1k,\]
and this shows that when the event 
$A^n_{\epsilon,p_1,p_2}$ holds, then $x$ gets aligned 
above or on $p_1x-p_1\epsilon n-p_1k$. 
This finishes proving that
the event $A^n_{\epsilon,p_1,p_2}$ is a sub-event of 
$D^n_a$.  Therefore, by \eqref{mainineq}, 
\[\mathbb{P}\left((D^{n}_a)^c\right)
\le \mathbb{P}\left((A^{n}_{\epsilon,p_1,p_2})^c\right)\le 
\exp\left(-n \left(-\frac{1+\ln k}{k}+\frac{\delta^2\epsilon^2}{16}\right)\right).\]
Similarly, and symmetrizing the above arguments,  
\[\mathbb{P}\left((D^{n}_b)^c\right)
\le \mathbb{P}\left((A^{n}_{\epsilon,p_1,p_2})^c\right)\le 
\exp\left(-n \left(-\frac{1+\ln k}{k}+\frac{\delta^2\epsilon^2}{16}\right)\right),\]
finishing, via \eqref{bbp}, the proof of the theorem.  
\hfill\end{proof}

\

Theorem~\ref{diagol} should prove useful in reducing the time to compute  
the LCS of two random sequences.  Indeed, first by \eqref{double}, when rescaled by $n$, 
the two-dimensional representation of an optimal alignment is, with high
probability and up to a distance of order $\epsilon>0$, above the line 
$x\rightarrow p_1x$ and below the line $x\rightarrow p_2x$. 
Moreover, $\epsilon>0$ can be taken as small
as we want, leaving it fixed though when $n$ goes to
infinity.  Next, simulations seem to indicate that 
the mean curve $\tilde\gamma$ is strictly concave at $p=1$.  
If strict concavity indeed hold, then $p_1$ and, say, $p_2= 1/p_1$, can 
be taken as close to $1$ as we want, and still satisfy the conditions of the 
theorem.  That is, taking $\epsilon$ as close to $0$ as we want 
and $p_1$ as close to $1$ as we want, the re-scaled 
two-dimensional representation of the optimal alignments would get 
uniformly as close to the diagonal as we want, as $n$ grows without bound.


\begin{figure}
\begin{center}
\includegraphics[width=5truein]{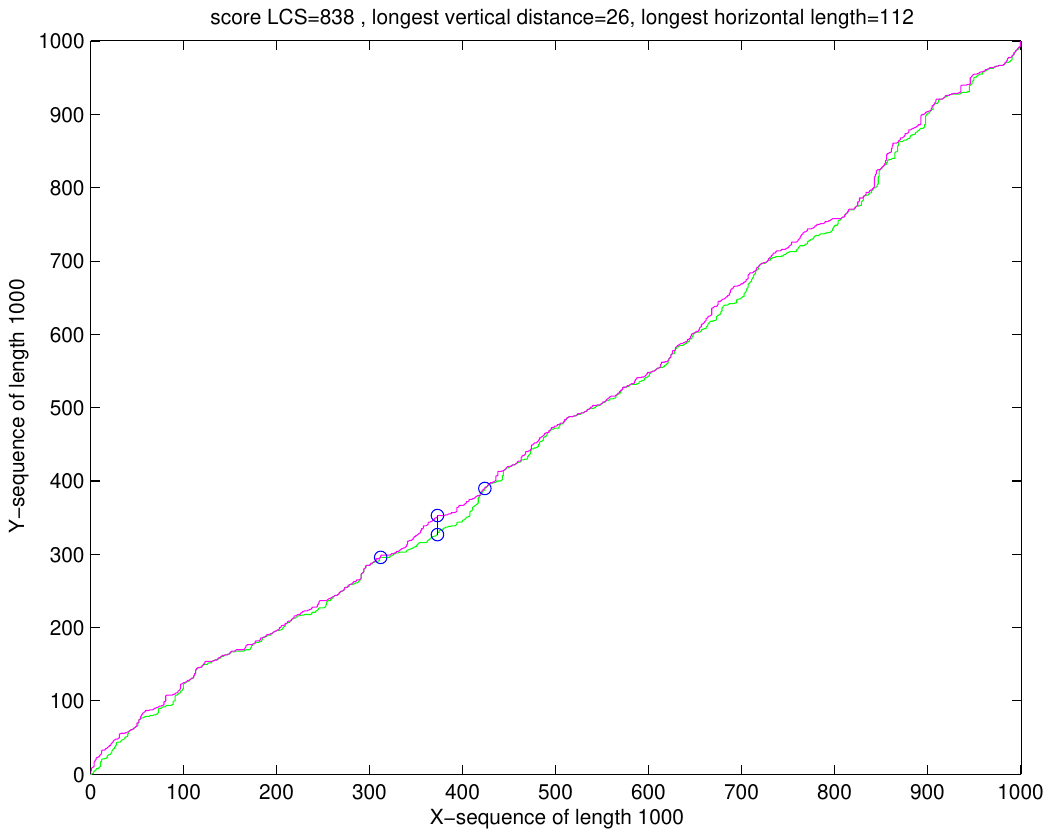}
\caption{$n=1000$, uniform Bernoulli sequences}
\end{center}
\label{fig:}
\end{figure}

Figure~1 is the graph of a simulation with two
iid binary sequences of length $n=1000$.  
All the optimal alignments are contained between the two 
graphs below and are thus seen to 
all stay extremely close to the diagonal.  
The maximal vertical distance between two optimal 
paths is $26$ and,  
for this vertical distance, 
the maximal horizontal stretch between which the two optimal paths 
split and then meet again is 
$112$.      

\section{Short string-lengths properties are generic}
\label{sectionproperty}

Often, a desirable property we want 
string-pairs to verify, e.g., a similar number of a given symbol or  
pattern, the presence of dominant matches, only holds with high probability  
and if the two strings have their lengths not too far away 
from each other.  Moreover, short strings are also often 
used as ``seeds" to find 
longer or more global similarities and homologous properties (\cite{SK} contains 
many examples of such instances in applied problems).  It is our purpose now to attempt to 
quantify such a generic phenomenon.  

To to so, let $\mathcal{P}$ be a relation 
assigning to every pair of strings $(x,y)$ 
the value $1$ if the pair $(x,y)$ has a certain property, and 
$0$ otherwise. Hence, if ${\mathcal A}$ is the alphabet 
we consider, 
\[\mathcal{P}: (\cup_k\mathcal{A}^k)\times (\cup_k \mathcal{A}^k)
\rightarrow \{0,1\},\] 
and if $\mathcal{P}(x,y)=1$, the string pair
$(x,y)$ is said to have the property $\mathcal{P}$.

Let now $\epsilon>0$, be fixed, and let $\vec{r}=(r_0,r_1,\dots,r_m)$ satisfy 
the condition \eqref{increasing}.  
Let also $B^n_{\mathcal{P}}(\vec{r},\epsilon)$ be the event 
that there is a proportion of at least $1-\epsilon$ 
of the string pairs 
\begin{equation} 
\label{stringpairs}
(X_{(i-1)k+1}\cdots X_{ik};Y_{r_{i-1}+1}\cdots Y_{r_i})
\end{equation}
satisfying the property $\mathcal{P}$, i.e., 
\begin{equation}\label{defB}
B^n_{\mathcal{P}}(\vec{r},\epsilon) = \left\{\sum_{i=1}^m 
\mathcal{P}(X_{(i-1)k+1}....X_{ik};Y_{r_{i-1}+1}\cdots Y_{r_i})\geq 
(1-\epsilon)m\right\}.
\end{equation}
Next, let $B^n_{\mathcal{P}}(\epsilon)$ be the event that for 
every optimal alignment 
the proportion of aligned string pairs \eqref{stringpairs} 
satisfying the property $\mathcal{P}$ is at least $1-\epsilon$, i.e., 
$B^n_{\mathcal{P}}(\epsilon)$
holds if and only if for every $\vec{r}=(r_0,r_1,\dots,r_m)$
satisfying \eqref{increasing} and such that 
$LC_n=L_n(\vec{r})$, the event 
$B^n_{\mathcal{P}}(\vec{r},\epsilon)$
holds.  
Finally, assume that as soon as $r_{i}-r_{i-1}\in [kp_1,kp_2]$,
the probability that the string-pairs
\eqref{stringpairs} have the required property
is at least $1-q$, $q\in [0, 1]$.  Hence, assume that
for every integer $\ell\in [kp_1,kp_2]$:
\[\bbp\left(\;\mathcal{P}(X_1\cdots X_k;Y_1Y_2\cdots Y_{\ell})=1\;\right)
\geq 1-q.\]
We investigate, now, how small 
$q=q(k)$ needs to be in order to ensure that a large proportion 
of the aligned string pairs \eqref{stringpairs} 
has the property $\mathcal{P}$ (for every optimal alignment). 
Recall that $A^n_{\epsilon,p_1,p_2}$ is the event 
that  every optimal alignment aligns a proportion of at least $1-\epsilon$ 
of the sub-strings $X_{(i-1)k+1}\cdots X_{ik}$ with sub-strings 
of $Y$ with length in $[kp_1,kp_2]$.  
Recall also that $\mathcal{R}_{\epsilon,p_1,p_2}$ 
is the set of integer vectors $\vec{r}=(r_0,r_1,\dots,r_m)$, 
satisfying \eqref{increasing} and such that there 
is at least $(1-\epsilon)m$ of the differences 
$r_i-r_{i-1}$ in $[kp_1,kp_2]$. 

Below, we deal with a small modification of the event
$B^n_{\mathcal{P}}(\vec{r},\epsilon)$. 
For this, let $\tilde{B}^n_{\mathcal{P}}(\vec{r},\epsilon)$
be the event that among the aligned string pieces 
\eqref{stringpairs} there are no more than $m\epsilon$ which do not satisfy
the property $\mathcal{P}$ and have their length 
$r_i-r_{i-1}$ in $ [kp_1,kp_2]$. Clearly, 
for $\epsilon_1 > 0$, $\epsilon_2 > 0$, 
\[A^n_{\epsilon_1,p_1,p_2}\cap 
\left(\bigcap_{\vec{r}\in\mathcal{R}_{\epsilon_1,p_1,p_2}}
 \tilde{B}^n_{\mathcal{P}}(\vec{r},\epsilon_2) \right)
\subset B^n_{\mathcal{P}}(\epsilon_1+\epsilon_2),\]
and so
\begin{equation}
\label{bound4B2}\bbp((B^{n}_{\mathcal{P}}(\epsilon_1+\epsilon_2))^c)\le
\bbp((A^{n}_{\epsilon_1,p_1,p_2})^c)+
\sum_{\vec{r}\in\mathcal{R}_{\epsilon_1,p_1,p_2}}
\bbp((\tilde{B}^{n}_{\mathcal{P}}(\vec{r},\epsilon_2))^c).  
\end{equation}
Next, 
\[\bbp((\tilde{B}^{n}_{\mathcal{P}}(\vec{r},\epsilon_2))^c)
\le
\binom{m}{\epsilon_2 m}q^{\epsilon_2 m}  \le
\exp(H_e(\epsilon_2)m)q^{\epsilon_2m},\] where $H_e$ is the base $e$ entropy 
function, given by $H_e(x) = -x\ln x - (1-x)\ln(1-x), 0 < x < 1$. Hence, 
\begin{equation}
\label{bound4B}
\bbp((\tilde{B}^{n}_{\mathcal{P}}(\vec{r},\epsilon_2))^c)
\le
q^{\epsilon_2 m}\exp(H_e(\epsilon_2) m).
\end{equation}
Using \eqref{bound4B} into
\eqref{bound4B2} and, proceeding as in \eqref{inebin}, noting that 
 ${\mathcal{R}}_{\epsilon,p_1,p_2}$ has at most $(ek)^m$ elements, lead to 
\begin{equation}
\label{bound4B3}\bbp((B^{n}_{\mathcal{P}}(\epsilon_1+\epsilon_2))^c)\le
\bbp((A^{n}_{\epsilon_1,p_1,p_2})^c)+
(ek)^m q^{\epsilon_2m}\exp(H_e(\epsilon_2)m).  
\end{equation}
Taking $q(k) = 1/(2ek)^{1/\epsilon_2}$, 
finally yields
\begin{equation}
\label{bound4B4}\bbp((B^{n}_{\mathcal{P}}(\epsilon_1+\epsilon_2))^c)\le
\bbp((A^{n}_{\epsilon_1,p_1,p_2})^c)+
\exp(\left(H_e(\epsilon_2)-\ln 2)m\right).
\end{equation}

\

But, for $\epsilon_2<1/2$, $H_e(\epsilon_2)<\ln 2$, and then  
$\exp(\left(H_e(\epsilon_2)-\ln 2)m\right)$
is exponentially small in $m$.  Now, our main theorem provides an exponentially small lower bound 
on $\bbp((A^{n}_{\epsilon_1,p_1,p_2})^c)$.   Therefore, \eqref{bound4B4}  
asserts that a high proportion of the aligned 
string pairs \eqref{stringpairs} has property $\mathcal{P}$, in any optimal 
alignment, as soon as for pairs \eqref{stringpairs} with similar length, 
$q(k) \le {1}/{(2ek)^{1/\epsilon_2}}$, where
\[q(k):=\max_{\ell\in[kp_1,kp_2]}
\bbp(\text{the pair }(X_1\cdots X_k;Y_1\cdots Y_{\ell})\text{ does not satisfy property }
\mathcal{P}).\]
These assertions are summarized in the next theorem, which is obtained 
by letting, above, $\epsilon_1=\epsilon_2=\epsilon/2$, using also Theorem~\ref{maintheorem}.    
\begin{theorem}
\label{lasttheorem}
Let $0< \epsilon < 1$.  Let $0< p_1<1<p_2$ be such that $\tilde\gamma(p_1)< \gamma^*$ and
$\tilde\gamma(p_2)<\gamma^*$, and let 
$0< \delta < \min\left(\gamma^*-\tilde\gamma(p_1), \gamma^*-\tilde\gamma(p_2)\right)$. 
Finally, let the integer $k\geq 1$ be such that 
\[\max_{\ell\in[kp_1,kp_2]}\bbp((X_1\cdots X_k;Y_1\cdots Y_\ell)
\text{ does not satisfy property $\mathcal{P}$})
\le \frac{1}{(2ek)^{2/\epsilon}}.\] 
Then, for any optimal alignment $\vec{r}$ (i.e., such that
$LC_n=L_n(\vec{r})$), the proportion 
of string pairs $(X_{(i-1)k+1}\cdots X_{ik};Y_{r_{i-1}+1}\cdots Y_{r_i})$
satisfying property $\mathcal{P}$ is at least $1-\epsilon$ with 
probability at least equal to:
\[
1-\bbp((A^{n}_{\epsilon/2,p_1,p_2})^c)-
\exp\left(\left(H_e\left(\frac{\epsilon}{2}\right)-\ln 2\right)m\right), 
\]
and thus at least equal to:
\[
1- \exp\left(-n\left(-\frac{1+\ln k}{k}+\frac{\delta^2\epsilon^2}{64}\right)\right)-
\exp\left(\frac{n}{k}\left(H_e\left(\frac{\epsilon}{2}\right)-\ln 2\right)\right),
\]
for all $n=n(\epsilon, \delta)$ large enough.  
\end{theorem}

\

Hence, from the above statement, the probability that less than 
a proportion $1-\epsilon$ of string pairs \eqref{stringpairs} 
have property $\mathcal{P}$ in every optimal alignment 
is exponentially small in $n$ (while holding 
$k$, $\epsilon$ and $\delta$ fixed) as soon as 
\begin{equation}
\label{conditionk}
k> \frac{64(1+\ln k)}{\epsilon^2\delta^2}, 
\end{equation}
and
\begin{equation}
\label{conditionprobability}\max_{\ell\in[kp_1,kp_2]}
\bbp((X_1\cdots X_k;Y_1\cdots Y_{\ell})\;\text{does not satisfy 
property $\mathcal{P}$})\le \frac{1}{(2ek)^{2/\epsilon}}.
\end{equation}

The above theorem is very useful for showing that when 
a property holds for aligned string pairs 
with similar lengths, say of order $k$, then 
the property typically holds in most parts of the 
optimal alignment. From our experience, 
for most properties one is interested in, such as  
the study of dominant matches 
in optimal alignments, 
when $p_1$ and $p_2$ are close to $1$, but 
fixed, then the probability that 
\[(X_1\cdots X_k;Y_1Y_2\cdots Y_\ell)\] 
does not satisfy this property 
is approximately the same for all $\ell\in[kp_1,kp_2]$. 
In other words, the behavior of the alignment 
of $X_1\cdots X_k$ with $Y_1\cdots Y_\ell$, 
does not depend much on  $\ell$ 
as soon as $\ell$ is close to $k$ and $k$ is fixed.  
From \eqref{conditionprobability}, what is needed there is a bound, on the left 
hand-side probability, smaller than any inverse polynomial-order in $k$.  
(At least to be able to take $\epsilon$ 
as close to $0$ as one wants to.)  If instead $\epsilon>0$ is chosen  
small but fixed, then an inverse polynomial bound 
with a very large exponent will do).  So, if this probability 
is, for example, of order $k^{-\ln k}$ or $e^{-k^\alpha}$ 
for some constant $\alpha>0$, the 
condition \eqref{conditionprobability} is satisfied by taking 
$k$ large enough.  Similarly, 
condition \eqref{conditionk} is always satisfied 
for $k$ large enough.

We could also envision using Monte Carlo 
simulation to find a bound 
for the probability on the left of \eqref{conditionprobability}.  
For that purpose, assume that $\epsilon = 0.2$ and take $\delta=0.1$.  
Then, by \eqref{conditionk}, $k$ must be at least $2 518 253$.    
The probability that strings of length approximately $k$ 
do not satisfy property $\mathcal{P}$ must be at most 
$(2ek)^{-10}\approx (13 690 642)^{-10}$, so a probability  
smaller than $10^{-70}$.  However, this is hardly feasible, 
indeed, to show that a probability is as small as $10^{-70}$, 
one would need to run an order of $10^{70}$ simulations.  

\

\centerline {\bf Further Improvements}

\

There are several ways to improve our various bounds. First,
we took as upper bound for cardinality of $\mathcal{R}_{\epsilon,p_1,p_2}$
the value $\binom{n}{m}$, which can be improved 
as follows: first note that if $\vec{r}=(r_0,r_1,\dots,r_m)\in
\mathcal{R}_{\epsilon,p_1,p_2}$, then at least $(1-\epsilon)m$
of the lengths $r_{i+1}-r_i$ are in the interval $[kp_1,kp_2]$.
To determine these lengths we have at most
\begin{equation}
\label{boundI}
((p_2-p_1)k)^{m}
\end{equation} choices. Then, there can be as many
as $\epsilon m$ of the lengths $r_{i+1}-r_i$, which are
not in $[kp_1,kp_2]$. Choosing those lengths is like 
choosing at most $\epsilon m$ points from a set of at most
$n$ elements. Hence, we get as upper bound $\binom{n}{\epsilon m}$
which, in turn, can be upper bounded by, say, 
\begin{equation}
\label{boundII}
\left(\frac{ek}{\epsilon} \right)^{\epsilon m}, 
\end{equation}
or via the entropy bound $\exp(nH_e(\epsilon/k))$.  
Finally, we have to decide which among the $m$ lengths $r_i-r_{i-1}$
have their length in $[kp_1,kp_2]$ and which have not.
That choice is further bounded via:
\begin{equation}
\label{boundIII}
\binom{m}{\epsilon m}\le \exp(H_e(\epsilon) m).
\end{equation}
Combining the bounds \eqref{boundI}, \eqref{boundII} and \eqref{boundII},
yields
\begin{equation}
\label{goodboundR}
{\rm Card}\left(\mathcal{R}_{\epsilon,p_1,p_2} \right)
\le 
\left((p_2-p_1)k\left(\frac{ek}{\epsilon} \right)^\epsilon \exp(H_e(\epsilon))
\right)^m
\end{equation}

With this better bounding for the cardinality
of $\mathcal{R}_{\epsilon,p_1,p_2}$, the inequality \eqref{bound4B3} becomes:
\begin{align}
\label{nextbound}
\bbp((B^{n}_{\mathcal{P}}(\epsilon_1+\epsilon_2))^c)&\le
\bbp((A^{n}_{\epsilon_1,p_1,p_2})^c)\nonumber\\
&\quad +
\left((p_2-p_1)k\left(\frac{ek}{\epsilon_1} \right)^{\epsilon_1} 
\exp(H_e(\epsilon_1))
\right)^m q^{\epsilon_2m}\exp(H_e(\epsilon_2)m), 
\end{align}
which when combined with Theorem~\ref{maintheorem} yields that
\begin{align}
\bbp((B^{n}_{\mathcal{P}}(\epsilon_1+\epsilon_2))^c)&\le
\label{nextbound2}
\exp\left(-n\left(-\frac{1 + \ln k}{k}+\frac{\delta^2\epsilon_1^2}{16}\right)\right)\nonumber \\
&\quad +
\left((p_2-p_1)k\left(\frac{ek}{\epsilon_1} \right)^{\epsilon_1} 
\exp(H_e(\epsilon_1))
 q^{\epsilon_2}\exp(H_e(\epsilon_2))\right)^m.
\end{align}

Again, this last expression is exponentially small 
in $n$ (assuming $k$ fixed) if the following two
conditions are satisfied:

\noindent
(i) \begin{equation}\label{k}
k>\frac{16(1+\ln k)}{\epsilon_1^2\delta^2},
\end{equation}
(the smallest integer $k$ satisfying \eqref{k} with $\epsilon_1 = 0.2$ and $\delta = 0.1$ is now $570 146$) 
and 
(ii) 
\begin{equation}\label{qk}
q(k)<\frac{1}{\left((p_2-p_1)k\left(\frac{ek}{\epsilon_1} \right)^{\epsilon_1} 
\exp(H_e(\epsilon_1)+H_e(\epsilon_2)) \right)^{1/\epsilon_2}};    
\end{equation}

\noindent
and combining these last two conditions yields:
\begin{equation}
\label{maincondition}
q(k)<
\left(
\frac{\epsilon_1^2\delta^2}{
(p_2-p_1)16(1+\ln k)\left(\frac{ek}{\epsilon_1} \right)^{\epsilon_1} 
\exp(H_e(\epsilon_1)+H_e(\epsilon_2))
}
\right)^{1/\epsilon_2}.
\end{equation}

Typically $\epsilon_1+\epsilon_2$ 
should be of a given order. So, let us maximize 
the right-hand side of \eqref{maincondition} under 
the constraint     
$\epsilon=\epsilon_1+\epsilon_2$.  To do so, note 
that the power $1/\epsilon_2$ has a much more minimizing influence 
than the expression $\epsilon_1^2$ in the numerator, while 
$1 \le \exp(H_e(\epsilon_1)+ H_e(\epsilon_2)) \le 2$ and so  
this last quantity does not have much of an influence.  Also, note 
$(ek/\epsilon_1)^{\epsilon_1}$ is somewhat negligible compared 
to $ek$. So, at first, let us disregard the quantities 
$(ek/\epsilon_1)^{\epsilon_1}$ and $\exp(H_e(\epsilon_1))+H_e(\epsilon_2))$,   
and let 
\[g(k,\epsilon_1,\epsilon_2):=
\left(
\frac{\epsilon_1^2\delta^2}{
(p_2-p_1)16(1+\ln k) 
}
\right)^{1/\epsilon_2}.
\]

Clearly, $g(k,\epsilon_1,\epsilon_2)$ is larger 
than the bound on the right-hand side of \eqref{maincondition} 
and when all the parameters $p_1, p_2$ and $\delta$ are held 
fixed, $g(k,\epsilon_1,\epsilon_2)$ is decreasing 
in both $\epsilon_1$ and $\epsilon_2$. However, 
$\epsilon_2$ ``has more decreasing influence'' than $\epsilon_1$.  
Therefore, given $\epsilon$ and given 
that all the parameters are fixed (including $k$), 
maximizing $g(k,\epsilon_1,\epsilon_2)$ 
under the constraint $\epsilon_1+\epsilon_2=\epsilon$, 
$\epsilon_1,\epsilon_2>0$ lead to a quantity 
where $\epsilon_2$ is quite a bit larger than $\epsilon_1$. 

Could Monte Carlo simulations be realistic with
$\epsilon=0.1$ and the bounds which we have? The answer is no.  
Indeed, at first, $\delta/(p_2-p_1)$ gets 
better when $|p_2-p_1|$ increases since the derivative of $\tilde\gamma$ 
at $p=1$ is zero.  When the interval $[kp_1,kp_2]$ 
becomes too large however, then the property might no longer 
hold with high probability for all pairs $(X_1\cdots X_k,Y_1\cdots Y_\ell)$, 
with $\ell\in[kp_1,kp_2]$.  So, we will take $[p_1, p_2]$, as large as possible, 
so this property still 
holds with high probability for all the string pairs 
mentioned before.  With such a choice, 
$\delta/(p_2-p_1)$ can be treated as a constant.  
Somewhat, optimistically, say that the constant is less than 
$1/3$.  Now if $\epsilon=0.1$, then $\epsilon_1,\epsilon_2\le 0.1$. 
In that case, 
\[g(k,\epsilon_1,\epsilon_2)\le g(k,0.1,0.1)\le
\left(
\frac{0.01\delta}{
3\times 16(1+\ln k) 
}
\right)^{10}.
\]
Returning to $\eqref{k}$ and taking $\delta=0.2$, 
we find that $k$ must be larger than $10^{10}$, 
so that $\ln k$ is bigger than $20$. 
With this in mind, and in the present case where 
$\epsilon_1+\epsilon_2=0.1$, we find that 
$g(k,\epsilon_1,\epsilon_2)$  
is smaller $10^{-56}$, so there is still 
little hope to perform Monte Carlo simulation here.  

Monte Carlo simulation with $\epsilon_1=0.1$ and $\epsilon_2=0.2$:  
Take also $\delta=0.2$ and $\delta/(p_2-p_1)=1/2$.  
With these values, and using \eqref{k}, then $k$ must be somewhat 
larger than ${10}^2\cdot 10 \cdot 8\ln (24000)\approx 10^5$. 
Then, by \eqref{qk}, $q(k)$  should also be less than 
\[\left(10^2\cdot 10\cdot 11\cdot 8  \right)^{-5}\approx 10^{-25}.\] 
This is still a difficult order for Monte Carlo simulation and 
if we had $\epsilon_2=0.3$ instead, then 
we would get a bound $10^{-15}$ which remains a difficult order.  

When only dealing with the inequality \eqref{qk}, things look somewhat better.
Take $k=1000$ and $(p_1-p_1)k=100$, then the bound on $q(k)$ 
is of order about $10^{-5}$ which is feasible with Monte Carlo.
So, if we could find another method than the one described 
here to make sure that most of the pieces of strings 
$X_{(i-1)k+1}X_{(i-1)k+2}\cdots X_{ik}$ are aligned with pieces of 
similar length we would end up in a favorable setting.

\

\noindent
{\bf Acknowledgments:} Many thanks 
to the referee for a thoughtful and detailed reading of this manuscript.

\end{document}